\documentstyle[fleqn,12pt]{article}
\begin{document}\pagenumbering{arabic}\setcounter{page}{1}
\pagestyle{plain}\baselineskip=16pt

\rightline{MSUMB-97-03, April 97}
\begin{center}
{\large\bf Two-parameter nonstandard deformation of 2x2 matrices} 
\end{center}

\vspace{1cm}\noindent
Salih Celik 

\noindent
{\footnotesize Mimar Sinan University, Department of Mathematics, 
80690 Besiktas, Istanbul, TURKEY. }

\noindent
Sultan A. Celik 

\noindent
{\footnotesize Yildiz Technical University, Department of Mathematics, 
Sisli, Istanbul, TURKEY.}

\vspace{2.5cm}
\noindent
{\bf Abstract}

\noindent
We introduce a two-parameter deformation of 2x2 matrices without 
imposing any condition on the matrices and give the universal R-matrix 
of the nonstandard quantum group which satisfies the quantum Yang-Baxter 
relation. Although in the standard two-parameter deformation the 
quantum determinant is not central, in the nonstandard case it is central. 
We note that the quantum group thus obtained is related 
to the quantum supergroup $GL_{p,q}(1|1)$ by a transformation.

\vfill\eject
\noindent
{\bf I. INTRODUCTION}

Recently the matrix groups of all 2x2 nonsingular matrices like $GL(2)$, 
$GL(1|1)$, {\it etc.}, were generalized in two ways as the standard 
deformation $^{1-3}$ and $h$-deformation.$^{4-7}$ Both are based on 
the deformation of the algebra of functions on the groups generated by 
coordinate functions that commute. 

In standard deformation of matrix groups, these commutation relations 
are determined by a matrix $R$ so that the functions do not commute 
but satisfy the equation 
$$\hat{R}(T \otimes T) = (T \otimes T) \hat{R},$$
such that, they coincide with the matrix groups for particular values 
of the deformation parameter. In the $h$-deformation, this property is 
the same as the standard deformation. 
The structure of the matrix groups is important in both deformations 
since the classical (or super) matrix groups are obtained in some limit 
of the deformation parameters. In this work we shall construct a two-parameter 
deformation of 2x2 matrices without imposing any such condition on the 
matrices just as in Ref. 8 and obtain a two-parameter generalization of their 
results. 

We briefly describe the content of this work. In section II we introduce 
the group $G_{p,q}$ of the 2x2 matrices by using an $R$ matrix. Section III 
is devoted to the corresponding Hopf algebra. In Sec. IV we give the universal 
enveloping algebra of this nonstandard quantum group. 

\noindent
{\bf II. G$_{p,q}$-MATRICES} 

Let 
$$T = \left(\matrix{ a & b \cr c & d \cr}\right) $$
be a 2x2 matrix with entries belonging to an algebra ${\cal A}$. We assume 
that the quantum group equation (no-grading) 
$$\hat{R} (T \otimes T) = (T \otimes T) \hat{R} \eqno(1)$$
holds, where 
$$\hat{R} = \left(\matrix{ 
-q &        0   & 0 & 0 \cr 
 0 & p^{-1} - q & 1 & 0 \cr 
 0 & qp^{-1}    & 0 & 0 \cr 
 0 &          0 & 0 & p^{-1} \cr} \right). \eqno(2)$$
Equation (1) explicitly gives the following relations: 
$$ a b = - q b a,  \quad d b = q b d, $$
$$ a c = - p c a,  \quad d c = p c d, \eqno(3) $$
$$ b c = pq^{-1} c b,  \quad  b^2 = 0 = c^2, $$
$$ ad =  da + (p^{-1} - q) b c, $$
where $p$ and $q$ are non-zero complex numbers with $pq \pm 1 \neq 0$. 

It can be checked that the matrix $R = P \hat{R}$, where $P$ is the usual 
permutation matrix, satisfies the quantum Yang-Baxter equation 
$$R_{12} R_{13} R_{23} = R_{23} R_{13} R_{12} \eqno(4)$$
and the matrix $\hat{R}$ satisfies the braid group equation 
$$\hat{R}_{12} \hat{R}_{23} \hat{R}_{12} = 
  \hat{R}_{23} \hat{R}_{12} \hat{R}_{23}. \eqno(5)$$

We now assume that the matrix elements $a$ and $d$ of $T$ are invertible. 
Then it is possible to define the inverse of $T$. To this end, we introduce 
$$ \Delta_1 = ad -p^{-1} b c, \quad 
   \Delta_2 = da + q^{-1} c b. \eqno(6) $$
Then one obtains 
$$ T^{-1}_R = \left(\matrix{ 
 \Delta_1^{-1} d    & - q \Delta_1^{-1} b \cr
  p \Delta_2^{-1} c &  \Delta_2^{-1} a \cr} \right),\eqno(7\mbox{a})$$
as the right inverse of $T$. After some calculations we get 
$$\Delta_1 d = d \Delta_1, \qquad \Delta_2 a = a \Delta_2,$$
$$\Delta_k b = - q^2 b \Delta_k, \qquad \Delta_k c = - p^2 c \Delta_k, \quad 
  k = 1, 2.$$
Using these relations we obtain 
$$T_L^{-1} = \left(\matrix{ 
  d\Delta_1^{-1}            & q^{-1} b \Delta_2^{-1} \cr
  - p^{-1} c  \Delta_1^{-1} & a \Delta_2^{-1} \cr} \right) = T^{-1}_R. 
\eqno(7\mbox{b}) $$
Thus the proper left and right inverses of $T$ are equal. 

It is easily verified that $a^2 \Delta^{-1}_2$ for all values of $p$ and $q$, 
commutes with $a$, $d$, and anticommutes with $b$, $c$. Furthermore 
$a^2 \Delta^{-1}_2$ is invertible. Therefore we obtain 
$$S(T) = T^{-1} = \left(\matrix{ 
 d^{-1}            & - a^{-1} b a^{-1} \cr
 - d^{-1} c d^{-1} &  a^{-1} \cr} \right) \left(\matrix{ 
  d^2 \Delta_1^{-1} &        0 \cr
   0                & a^2 \Delta_2^{-1} \cr} \right). \eqno(8) $$

We now consider the element 
$$ D(T) = ad^{-1} - b d^{-1} c d^{-1} = a^2 \Delta^{-1}_2. \eqno(9)$$
$D(T)$ cannot be regarded as a quantum determinant since it anticommutes 
with $b$ and $c$. However, we may regard the element 
$${\bf D}(T) = a [D(T) - d^{-1} b d^{-1} c] d^{-1} = [D(T)]^2 \eqno(10)$$
as the quantum determinant of $T$ where $D(T)$ is given by (9). 

In fact, it is easy to check that the matrix elements of the product 
matrix $T T'$ satisfy relations (3) for any two commuting quantum 
matrices $T$ and $T'$ whose elements obey (3). As a consequence of 
this argument, we have the following relation: 
$$D(T T') = D(T) D(T'). $$
This result means that ${\bf D}(T)$ is central. 

This case appears strange from the point of view of quantum group theory.$^2$ 
However it becomes clear from the point of view of the corresponding 
two-parameter quantum supergroup.$^9$ We know, from the work of Ref. 9, that 
the quantum superdeterminant of any supermatrix in GL$_{p,q}(1\vert 1)$ 
belongs to the centre of the algebra generated by the matrix elements of the 
supermatrix. In the Appendix, we shall show that this nonstandard quantum 
group is related to the quantum supergroup GL$_{p,q}(1\vert 1)$ by a 
transformation. So we may expect that the quantum superdeterminant in 
two-parameter nonstandard deformation must again be a central element. 

Now let the $n$-th power of $T$ be 
$$ T^n = \left ( \matrix{ A_n & B_n \cr
                          C_n & D_n \cr}
\right). \eqno(11) $$
Then it is easy to check the following relations: 
$$ A_n B_n = - q^n B_n A_n, \quad D_n B_n = q^n B_n D_n, $$
$$ A_n C_n = - p^n C_n A_n, \quad D_n C_n = p^n C_n D_n, \eqno(12)$$
$$ B_n^2 = 0 = C_n^2, \quad q^n B_n C_n = p^n C_n B_n,$$
and 
$$ A_n D_n = D_n A_n + (p^{-n} - q^n) C_n B_n. \eqno(13)$$
The proof of relation (13) is rather lengthy but straightforward. 

Let us finally note the following. If the sum $T + T'$ of two 
$G_{p,q}$ matrices $T$ and $T'$ is required to be a $G_{p,q}$ matrix 
then the equation 
$$\hat{R}' (T \otimes T') = (T' \otimes T) \hat{R}^{-1} \eqno(14)$$
holds, where 
$$\hat{R}' = \hat{R}^{-1} - (p - q^{-1}) I. \eqno(15)$$ 
Equation (14) explicitly reads 
$$a' a = pq a a', \quad d d' = pq d' d, \quad d' a = a d', $$
$$b' a = - p a b', \quad c' a = - q a c', \quad b d' = p d' b, \eqno(16)$$
$$b b' = - b' b, \quad c c' = - c' c \quad cd' = q d' c, $$
$$ a' b = - q b a' + (p^{-1} - q) b' a, \quad 
   a' c = - p c a' + (pq - 1) a c', $$
$$ b' c = pq^{-1} c b' + (q^{-1} - p) a d', \quad 
   b c' = pq^{-1} c' b + (p - q^{-1}) d' a, $$
$$ d b' = q b' d + (q - p^{-1}) b d', \quad 
   d c' = p c' d + (pq - 1) d' c, $$
$$a' d = d a' + (p^{-1} - q) (b c' + b' c).   $$
Note that the matrix $R' = P \hat{R}'$ again satisfies the quantum Yang-Baxter 
relation (4), where $P$ is the usual permutation matrix. 

\noindent
{\bf III. THE HOPF ALGEBRA STRUCTURE OF G$_{p,q}$ }

Let ${\cal A}$ be an algebra generated by the elements $a$, $b$, $c$ and $d$ 
satisfying the relations (3). Then ${\cal A}$ is the quoitent algebra 
$${\cal A} = {\cal C}[a,b,c,d]/J, $$
where ${\cal C}[a,b,c,d]$ is the free non-commutative algebra generated by 
$a$, $b$, $c$ and $d$ and $J$ is the ideal in ${\cal C}[a,b,c,d]$ 
generated by the relations (3). 

The usual co-product on the algebra ${\cal A}$ is defined by 
$$\Delta : {\cal A} \longrightarrow {\cal A} \otimes {\cal A}$$
such that 
$$\Delta(t^i{}_j) = t^i{}_k \otimes t^k{}_j, \quad T = (t^i_j) \eqno(17)$$
(sum over repeated indices) and the counit 
$$\varepsilon : {\cal A} \longrightarrow {\cal C} $$
such that 
$$\varepsilon(t^i{}_j) = \delta^i{}_j. \eqno(18)$$
The algebra ${\cal A}$ is now the matrix bialgebra generated by 1 and 
$T = (t^i_j)$, and it is a Hopf algebra with the antipode $S(T)$ 
which is given by (8). To give a proof of this, one has to verify the 
following: 
$$ (\mbox{id} \otimes \Delta) \circ \Delta = 
   (\Delta \otimes \mbox{id}) \circ \Delta, $$
$$ (\mbox{id} \otimes \varepsilon) \circ \Delta = 
   (\varepsilon \otimes \mbox{id}) \circ \Delta, \eqno(19)$$
$$ m \circ [(\mbox{id} \otimes S) \circ \Delta] = 
   m \circ [(S \otimes \mbox{id}) \circ \Delta], $$
where $m$ denotes the multiplication mapping 
$$m(a \otimes b) = a b $$
for any $a, b \in {\cal A}$. The proof follows directly. 

\noindent
{\bf IV. UNIVERSAL ENVELOPING ALGEBRA OF G$_{p,q}$} 

In this section we shall construct the quantum enveloping algebra in 
analogy with the FRT approach.$^2$ 

We consider the matrices $L^{\pm}$ with the generators $U_{\pm}$, 
$V_{\pm}$ and $X_{\pm}$, 
$$L^+ = \left(\matrix{ 
 U_+ & \lambda X_+ \cr 
 0   &  V_+ \cr} \right), \quad 
L^- = \left(\matrix{ 
 U_-           & 0  \cr 
 - \lambda X_- &  V_- \cr} \right), \eqno(20)$$
where $\lambda = q - p^{-1}$. The matrices $L^{\pm}$ satisfy the following 
relations: 
$$ \hat{R} L^{\pm}_1 L^{\pm}_2 = L^{\pm}_2 L^{\pm}_1 \hat{R}, \eqno(21)$$
$$ \hat{R} L^+_1 L^-_2 = L^-_2 L^+_1 \hat{R}, \eqno(22)$$
where $L_1 = L \otimes I$ and $L_2 = I \otimes L$. These relations give 
$$ [U_+,U_-] = [V_+,V_-] = [U_{\pm},V_{\pm}] = 0, $$
$$ U_+ X_{\pm} = - q^{\mp 1} X_{\pm} U_+, \quad 
   V_+ X_{\pm} = q^{\mp 1} X_{\pm} V_+, $$
$$ U_- X_{\pm} = - p^{\pm 1} X_{\pm} U_-, \quad 
   V_- X_{\pm} = p^{\pm 1} X_{\pm} V_-, \eqno(23) $$
$$ X_+ X_- - qp^{-1} X_- X_+ = {{U_+ V_- - V_+ U_-}\over {q - p^{-1}}}, \quad 
   X_{\pm}^2 = 0.$$
The coproduct of the generators is given by 
$$\Delta(L^{\pm}) = L^{\pm} \dot{\otimes} L^{\pm} \eqno(24)$$
where $\dot{\otimes}$ denotes tensor product and matrix multiplication. 
Explicity, the action of the coproduct $\Delta$ on the generators is 
$$ \Delta(U_{\pm}) = U_{\pm} \otimes U_{\pm}, $$
$$ \Delta(V_{\pm}) = V_{\pm} \otimes V_{\pm}, $$
$$ \Delta(X_+) = X_+ \otimes U_+ + V_+ \otimes X_+, \eqno(25)$$
$$ \Delta(X_-) = X_- \otimes V_- + U_- \otimes X_-. $$
The counit is given by 
$$ \varepsilon(L^{\pm}) = I. \eqno(26)$$
Explicitly, 
$$ \varepsilon(U_{\pm}) = \varepsilon(V_{\pm}) = 1, $$
$$ \varepsilon(X_{\pm}) = 0. \eqno(27)$$
The coinverse is given by 
$$ S(U_{\pm}) = U_{\pm}^{-1}, \quad S(V_{\pm}) = V_{\pm}^{-1}, $$
$$ S(X_+) = - U_+^{-1} X_+ V_+^{-1}, \eqno(28) $$
$$ S(X_-) = V_-^{-1} X_- U_-^{-1}. $$
Therefore one can easily verify that the algebra 
${\cal U}_{p,q}(U_{\pm}, V_{\pm}, X_{\pm})$ is a Hopf algebra generated by 
1, $U_{\pm}$, $V_{\pm}$, $X_{\pm}$ satisfying the relations (23). 

The coproduct of $U_{\pm}$ and $V_{\pm}$ together with the fact that they 
commute implies that they can be written as exponentials of commuting 
operators, 
$$U_+ = q^{- {H\over 2}} p^{N\over 2}, \quad 
  U_- = p^{{H\over 2}} q^{- {N\over 2}}, $$
$$V_+ = q^{- {H\over 2}} p^{- {N\over 2}}, \quad 
  V_- = p^{H\over 2} q^{N\over 2}, \eqno(29)$$
$$[H, N] = 0.$$
The commutation relations of $U_{\pm}$ and $V_{\pm}$ with $X_{\pm}$ in terms 
of new generators give the following: 
$$[H, X_{\pm}] = \pm 2 X_{\pm}, \quad [N, X_{\pm}] = 0, $$
$$X_+ X_- - qp^{-1} X_- X_+ = \left({p\over q}\right)^{{H + 1}\over 2} [N]_{pq}
  \eqno(30)$$ 
where 
$$[N]_{pq} = {{(pq)^{N\over 2} - (pq)^{- N\over 2}}\over 
             {(pq)^{1\over 2} - (pq)^{- 1\over 2}}}. \eqno(31)$$
Moreover, the coproduct is now 
$$\Delta(H) = H \otimes {\bf 1} + {\bf 1} \otimes H, $$
$$\Delta(N) = N \otimes {\bf 1} + {\bf 1} \otimes N, \eqno(32)$$
$$\Delta(X_+) = X_+ \otimes q^{- {H\over 2}} p^{N\over 2} + 
                q^{- {H\over 2}} p^{- {N\over 2}} \otimes X_+, $$
$$\Delta(X_-) = X_- \otimes p^{H\over 2} q^{N\over 2} + 
                p^{H\over 2} q^{- {N\over 2}} \otimes X_-. $$

\vfill\eject
\noindent
{\bf ACKNOWLEDGMENTS}

This work was supported in part by T. B. T. A. K. the Turkish Scientific 
and Technical Council. We would like to express our deep gratitude to the 
referee for critical comments on the manuscript. 

\noindent
{\bf V. APPENDIX} 

\noindent
{\bf A. Nonstandard quantum planes} 

In this section, we shall consider quantum planes which are 
similar to quantum superplanes introduced by Manin.$^3$ 

(1) Quantum plane $A_p$: This plane, or, rather the polynomial 
function ring on it is generated by coordinates $x$ and $\theta$ with 
the commutation rules 
$$x \theta = - p \theta x, \quad \theta^2 = 0, \eqno(A1)$$
where $p$ is a complex number. The coordinates anticommute for $p = 1$ and 
commute for $p = - 1$. 

(2) Quantum plane $A_q^*$: This plane is generated by 
coordinates $\varphi$ and $y$ with commutation rules 
$$ \varphi^2 = 0, \quad \varphi y = q^{-1} y \varphi \eqno(A2)$$
where $q$ is a complex number. The quantum plane $A^*_q$ is dual to 
the quantum plane $A_p$. 

Note that the relations (A1) and (A2) are equivalent to the relations 
$$\hat{R} (X \otimes X) = - q (X \otimes X), \quad 
  \hat{R} (Y \otimes Y) = p^{-1} (Y \otimes Y). $$

\noindent
{\bf B. Nonstandard quantum deformation of 2x2 matrices with nonstandard 
quantum planes} 

Let $G$ be a matrix Lie group of rank 2 and $T$ be any element of $G$, i.e., 
$$T = \left(\matrix{ 
  a & b \cr c & d \cr}\right)$$
with entries belonging to an algebra ${\cal A}$. 

We consider linear transformations $T$ with the following properties: 
$$T: A_p \longrightarrow A_p, \quad 
  T: A_q^* \longrightarrow A_q^*. \eqno(A3)$$
The action of $T$ on points of $A_p$ and $A_q^*$ is 
$$\left(\matrix{\overline{x} \cr \overline{\theta} \cr}\right) = 
  \left(\matrix{ 
  a & b \cr c & d \cr}\right)
\left(\matrix{ x \cr \theta \cr}\right), \quad 
\left(\matrix{\overline{\varphi} \cr \overline{y} \cr}\right) = 
  \left(\matrix{ 
  a & b \cr c & d \cr}\right)
\left(\matrix{ \varphi \cr y \cr}\right). \eqno(A4)$$
We assume that the matrix elements of $T$ commute with the coordinates of 
$A_p$ and $A_q^*$. As a consequence of the linear transformations in (A3) 
the vectors $\left(\matrix{\overline{x} \cr \overline{\theta} \cr}\right)$ and 
$\left(\matrix{\overline{\varphi} \cr \overline{y} \cr}\right)$ should belong 
to $A_p$ and $A_q^*$, respectively. This imposes $(p,q)$-commutation 
relations among the entries of $T$ in (3). 

Note that it can be checked that the maps 
$$\delta: A_p \longrightarrow G \otimes A_p, \quad 
  \delta^*: A_q^* \longrightarrow G \otimes A_q^* \eqno(A5)$$
such that 
$$\delta(X) = T \otimes X, \quad \mbox{i.e.} \quad \
  delta(x_i) = t^j_i \otimes x_j, \quad 
  X = \left(\matrix{ x \cr \theta \cr}\right) $$
$$\delta(Y) = T \otimes Y, \quad \mbox{i.e.} \quad 
  \delta(y_i) = t^j_i \otimes y_j, \quad 
  Y = \left(\matrix{ \varphi \cr y \cr}\right) \eqno(A6)$$
define the co-action of the quantum group $G_{p,q}$ on 
the nonstandard quantum planes $A_p$ and $A_q^*$, respectively. 

Finally, one can show that the matrix quantum group (3) is isomorphic 
to the quantum supergroup $GL_{p,q}(1|1)$. Indeed, if we define the 
transformation 
$$T' = T D \eqno(A7)$$
where $T$ is a matrix whose the matrix elements satisfy (3) and 
$$D = \left(\matrix{ 1 & 0 \cr 0 & g \cr} \right), \quad D^2 = I \eqno(A8)$$
and we assume that $g$ commutes with $a$ and $d$, and anticommutes with $b$ 
and $c$, then $T' \in GL_{p,q}(1|1)$ as discussed in Ref. 9. 
In this case 
$\Delta(T') = T' \dot{\otimes} T'$, {\it etc.}, are unchanged. 
One easily sees that when $p = q$, these relations go back to those of Ref. 8. 

{\footnotesize 
\noindent
$^1$ V. G. Drinfeld, Proceedings of the IMS, Berkeley, 1986 (unpublished). \\
$^2$ N. Y. Reshetikhin, L. A. Takhtajan and L. D. Faddeev, 
     Leningrad Math. J. {\bf 1}, 193 (1990).\\
$^3$ Yu I. Manin, 
     Commun. Math. Phys. {\bf 123}, 163 (1989).\\ 
$^4$ S. L. Woronowicz, 
     Rep. Math. Phys. {\bf 30}, 259 (1991). \\
$^5$ S. Zakrzewski, 
     Lett. Math. Phys. {\bf 22}, 287 (1991); \\
\hspace*{0.2cm} C. H. Ohn, 
     Lett. Math. Phys. {\bf 25}, 85 (1992); \\
\hspace*{0.2cm} B. A. Kupershmidt, 
     J. Phys. A: Math. Gen. {\bf 25}, L1239 (1992); \\
\hspace*{0.2cm} V. Karimipour, 
     Lett. Math. Phys. {\bf 30}, 87 (1994). \\
$^6$ A. Aghamohammadi, M. Khorrami and A. Shariati, 
     J. Phys. A: Math. Gen. {\bf 28}, L225 (1995).\\
$^7$ S. \c Celik, Lett. Math. Phys. {\bf 42}, 299 (1997), 
     and their references.  \\
$^8$ S. Majid and M. J. Rodriguez-Plaza, 
     J. Math. Phys. {\bf 36}, 7081 (1995). \\
$^{9}$ L. Dabrowski and L. Wang, 
       Phys. Lett. B {\bf 266}, 51 (1991). 
}

\end{document}